\documentclass[12pt]{article}
\usepackage{amssymb}

\input tcilatex

\begin{document}

\title{Hodge equations with change of type}
\author{Thomas H. Otway\thanks{%
email: otway@ymail.yu.edu} \\
\\
\textit{Departments of Mathematics and Physics,}\\
\textit{Yeshiva University, 500 W 185th Street,}\\
\textit{New York, New York 10033, USA}}
\date{}
\maketitle

\begin{abstract}
A geometric interpretation is given for certain elliptic-hyperbolic systems
in the plane. Among several examples, one which reduces in the elliptic
region to the equations for harmonic 1-forms on the projective disc is
studied in detail. \ A boundary-value problem for this example is formulated
and shown to possess weak solutions.\ \textit{MSC2000}: 35M10, 58J99. \ 
\end{abstract}

\section{Introduction}

Harmonic forms $u$ on a Riemannian manifold satisfy the \textit{Hodge
equations} 
\begin{equation}
\delta u=du=0,
\end{equation}
where $d$ is the exterior derivative and $\delta $ its adjoint. \ In the
case of 1-forms, these equations have the local form 
\begin{equation}
\left| G\right| ^{-1/2}\partial _{i}\left( G^{ij}\sqrt{\left| G\right| }%
u_{j}\right) =0,
\end{equation}
\begin{equation}
\partial _{i}u_{j}dx^{i}\wedge dx^{j}=\frac{1}{2}\left( \partial
_{i}u_{j}-\partial _{j}u_{i}\right) dx^{i}\wedge dx^{j}=0,
\end{equation}
where $G_{ij}$ is the metric tensor on the manifold. \ If\ eqs. (2), (3) are
defined on a singular 2-manifold, it may happen that the equations can be
rewritten as a system of mixed type in $\Bbb{R}^{2}$ in which the parabolic
curve lies along the singularity of the manifold. \ This yields a geometric
interpretation of certain elliptic-hyperbolic systems in the plane.

Perhaps the simplest example is a metric which changes from Euclidean to
Minkowskian along the $x$-axis. In this case the system (2), (3) reduces to
a potential equation on one side of the metric singularity and to a wave
equation on the other side, leading to a first-order system of the form 
\[
u_{1x}+sgn(y)u_{2y}=0, 
\]
\[
u_{1y}-u_{2x}=0. 
\]
This corresponds in the case $u_{1}=u_{x},$ $u_{2}=u_{y}$ to the
Lavrent'ev-Bitsadze equation.

An example possessing more interesting geometry can be constructed by taking
Beltrami's hyperbolic model [1] for the projective disc $\Bbb{P}^{2}$ as the
underlying surface. \ The metric tensor in this model is the matrix 
\[
G_{ij}=\frac{1}{\left( 1-x^{2}-y^{2}\right) ^{2}}\left[ 
\begin{array}{cc}
1-y^{2} & xy \\ 
xy & 1-x^{2}
\end{array}
\right] . 
\]
The matrix 
\[
G^{ij}=\left( 1-x^{2}-y^{2}\right) \left[ 
\begin{array}{cc}
1-x^{2} & -xy \\ 
-xy & 1-y^{2}
\end{array}
\right] 
\]
becomes indefinite, and the determinant 
\[
G=\frac{1}{1-x^{2}-y^{2}} 
\]
becomes singular, on the line at infinity of the model, which corresponds to
the circle $x^{2}+y^{2}=1$.

But the equations can be redefined so that the metric singularity on the
unit circle in $\Bbb{P}^{2}$ is replaced by a change of type on the unit
circle in $\Bbb{R}^{2}$. \ Writing out eq. (2) in coordinates, we obtain 
\[
\left( 1-x^{2}-y^{2}\right) \{\left[ \left( 1-x^{2}\right) u_{1}\right]
_{x}-\left( xyu_{1}\right) _{y}-\left( xyu_{2}\right) _{x} 
\]
\begin{equation}
+\left[ \left( 1-y^{2}\right) u_{2}\right] _{y}-\left( xu_{1}+yu_{2}\right)
\}=0.
\end{equation}
Equation (3) implies that 
\begin{equation}
\left( xyu_{1}\right) _{y}+\left( xyu_{2}\right)
_{x}=2xyu_{1y}+xu_{1}+yu_{2}.
\end{equation}

Outside the unit circle the projective disc model no longer applies, but
eqs. (4), (5) are well defined and possess wave-like solutions in which
disturbances propagate along null geodesics of \ the distance element 
\[
ds^{2}=\frac{\left( 1-y^{2}\right) dx^{2}+2xydxdy+\left( 1-x^{2}\right)
dy^{2}}{\left( 1-x^{2}-y^{2}\right) ^{2}}. 
\]
Borrowing the terminology of fluid dynamics, we call an expression such as $%
ds^{2}$ the \textit{flow metric} associated to a system such as (4), (5).

In order for a 1-form $u$ to satisfy (4), (5), it is sufficient for $u$ to
satisfy a system of first-order equations on $\Bbb{R}^{2}$ having the form 
\begin{equation}
Lu=g,
\end{equation}
where 
\[
L=\left( L_{1},L_{2}\right) ,\;g=\left( g_{1},g_{2}\right) , 
\]
\[
u=\left( u_{1}\left( x,y\right) ,u_{2}\left( x,y\right) \right) ,\;\left(
x,y\right) \in \Omega \subset \subset \Bbb{R}^{2}, 
\]
\begin{equation}
\left( Lu\right) _{1}=\left[ \left( 1-x^{2}\right) u_{1}\right]
_{x}-2xyu_{1y}+\left[ \left( 1-y^{2}\right) u_{2}\right]
_{y}-2xu_{1}-2yu_{2},
\end{equation}
and 
\[
\left( Lu\right) _{2}=u_{1y}-u_{2x}. 
\]
If $y^{2}\neq 1,$ we can replace the second component of $L$ by the
expression 
\begin{equation}
\left( Lu\right) _{2}=\left( 1-y^{2}\right) \left( u_{1y}-u_{2x}\right) ,
\end{equation}
which has the same annihilator.

The second-order terms of eqs. (6)-(8) can be written in the form $%
Au_{x}+Bu_{y},$ where 
\[
A=\left[ 
\begin{array}{cc}
1-x^{2} & 0 \\ 
0 & -\left( 1-y^{2}\right)
\end{array}
\right] 
\]
and 
\[
B=\left[ 
\begin{array}{cc}
-2xy & 1-y^{2} \\ 
1-y^{2} & 0
\end{array}
\right] . 
\]
If $y^{2}\neq 1$, the characteristic equation 
\[
\left| A-\lambda B\right| =-\left( 1-y^{2}\right) \left[ \left(
1-y^{2}\right) \lambda ^{2}+2xy\lambda +\left( 1-x^{2}\right) \right] 
\]
possesses two real roots $\lambda _{1},\lambda _{2}$ on $\Omega $ precisely
when $x^{2}+y^{2}>1$. \ Thus the system is elliptic in the intersection of $%
\Omega $ with the open unit disc centered at $\left( 0,0\right) $ and
hyperbolic in the intersection of $\Omega $ with the complement of the
closure of this disc. The boundary of the unit disc, along which this change
in type occurs, is the line at infinity in $\Bbb{P}^{2}$ and a line
singularity of the tensor $G_{ij}$.

L. K. Hua used variable separation, Poisson kernel, and D'Alembert methods
to solve boundary-value problems for a scalar equation which resembles the
system (6)-(8) [3]. \ Precisely, the scalar equation studied by Hua consists
of the conserved quantities in an equation which can be obtained from
(6)-(8) by choosing $u_{1}=u_{x}$, $u_{2}=u_{y}$, and $g_{1}=g_{2}=0$. \ A
form of the equation studied by Hua with $g_{1}\neq 0$ was solved by Ji and
Chen [4]. \ Inside the unit disc, these choices correspond to replacing the
Hodge operator on 1-forms with the Laplace-Beltrami operator on scalars,
modulo lower-order terms. \ We emphasize that eqs. (6)-(8), even without the
lower-order terms, are not equivalent to an equation of the form studied by
Hua if $g_{2}\neq 0$ or if the vector $\left( u_{1},u_{2}\right) $ is not
continuously differentiable. \ Beyond this, the form of the lower-order
terms which do not appear in [3] affect our analysis of the equations in
Secs. 3-6. \ (In Sec. 6 we consider the equations in the absence of
lower-order terms.) \ Finally, in Refs. 3 and 4 conditions are placed on
characteristics, as in the classical Tricomi problem; in Secs. 3-6
conditions are placed only on the noncharacteristic part of the boundary, as
in the classical Frankl' problem.

\section{Other systems of mixed type}

The analogy between the two systems (2), (3) and (6)-(8) can be extended to
other equations of mixed type, although generally these systems will have
less interesting geometry than the projective disc. \ For example, the
system introduced by Morawetz [5] as a vehicle for studying the Chaplygin
equations is of a broadly similar form, as is the system studied in Ref. 8.

\subsection{Equations of fluid dynamics}

The geometry of eqs. (6)-(8) is in some sense dual to that of a well known
transform of the velocity potential for transonic flow in the hodograph
plane. \ Denote by $\left( u_{1}\left( x,y\right) ,u_{2}\left( x,y\right)
\right) $ the velocity components of a steady flow expressed in coordinates $%
\left( x,y\right) $. \ The hodograph transformation introduces $u_{1}$, $%
u_{2}$ as independent coordinates. \ The continuity equations for the
velocity potential under standard simplifying assumptions can now be written
in the linear form ([2], eq. (3.6)) 
\[
\left( c^{2}-u_{1}^{2}\right) y_{u_{2}}+u_{1}u_{2}\left[ x_{u_{2}}+y_{u_{1}}%
\right] +\left( c^{2}-u_{2}^{2}\right) x_{u_{1}}=0, 
\]
\[
x_{u_{2}}-y_{u_{1}}=0. 
\]
Here 
\[
c^{2}=1-\frac{\gamma _{a}-1}{2}\left( u_{1}^{2}+u_{2}^{2}\right) , 
\]
where $\gamma _{a}>1$ is the adiabatic constant of the medium. \ This system
corresponds to eqs. (2), (3) where the parabolic curve is a circle of radius 
$\sqrt{2/\left( \gamma _{a}+1\right) }$ centered at the point $u_{1}=0,$ $%
u_{2}=0$ and the metric tensor in eq. (2) is the matrix 
\[
\widetilde{G}_{ij}=\frac{1}{c^{2}\left( c^{2}-u_{2}^{2}-u_{1}^{2}\right) }%
\left[ 
\begin{array}{cc}
c^{2}-u_{1}^{2} & -u_{1}u_{2} \\ 
-u_{1}u_{2} & c^{2}-u_{2}^{2}
\end{array}
\right] . 
\]
Consider for simplicity the lower limit of the range of values for $\gamma
_{a},$ in which $c^{2}$ is approximately normalized. \ In this artificially
simple case, the change of type occurs on the boundary of the unit circle
and the continuity equations in the hodograph plane reduce to a replacement
of the metric tensor $G_{ij}\left( u_{1},u_{2}\right) $ of eq. (2) by the
tensor $\left( 1-u_{1}^{2}-u_{2}^{2}\right) ^{-2}G^{_{ij}}\left(
u_{1},u_{2}\right) $ (ignoring lower-order terms). \ We obtain a
second-order scalar equation if we introduce a function $\chi \left(
u_{1},u_{2}\right) $ satisfying 
\[
x=\chi _{u_{1}},\;y=\chi _{u_{2}} 
\]
(\textit{c.f.} eq. (3.8) of Ref. 2). \ The characteristic curves of the
resulting equation are relatively complicated, as they are given by a family
of epicycloids which intersect the parabolic curve in a family of cusps. \
This leads to complicated boundary-value problems for the equation. \ By
contrast, the characteristic curves corresponding to the ``dual'' system
(6)-(8) are exceedingly simple, as they are given by the set of all tangent
lines to the unit disc. \ This leads in our case to relatively simple
boundary-value problems. \ How much can be said \textit{a priori} about
relations between solutions of the two sets of boundary-value problems is
not immediately clear, however.

Without its lower-order terms and after a trivial relabelling of
coordinates, the system (6)-(8) can be interpreted as the hodograph image of
a quasilinear system having the form 
\begin{equation}
\left( 1-u_{1}^{2}-u_{2}^{2}\right) ^{m}\left[ \left( 1-u_{2}^{2}\right)
u_{1x}+u_{1}u_{2}\left( u_{1y}+u_{2x}\right) +\left( 1-u_{1}^{2}\right)
u_{2y}\right] =0,
\end{equation}
\begin{equation}
u_{2x}-u_{1y}=0,
\end{equation}
for $m\in \Bbb{R}.$ \ If the components $u_{1}(x,y)$ and $u_{2}(x,y)$ are
continuously differentiable in $x$ and $y$, then there is a potential
function $\varphi \left( x,y\right) $ such that 
\[
d\varphi \left( x,y\right) =\varphi _{x}dx+\varphi _{y}dy=u_{1}dx+u_{2}dy 
\]
on any domain having trivial de Rham cohomology. \ If $m=-3/2$, then the
resulting equation is the \textit{Hodge dual} of the minimal surface
equation, in the sense of [9], eqs. (2.23)-(2.29). \ If $\left(
1-u_{1}^{2}-u_{2}^{2}\right) ^{m}\neq 0,$ then the flow metric for eqs. (9),
(10) is conformally equivalent to the metric 
\[
ds^{2}=dx^{2}+dy^{2}-\left( d\varphi \right) ^{2}. 
\]
By comparison, the flow metric for the gas dynamics equation 
\[
\left( 1-\frac{u_{1}^{2}}{c^{2}}\right) u_{1x}-\frac{u_{1}u_{2}}{c^{2}}%
\left( u_{1y}+u_{2x}\right) +\left( 1-\frac{u_{2}^{2}}{c^{2}}\right)
u_{2y}=0 
\]
is conformally equivalent to the metric 
\[
ds^{\prime 2}=dx^{2}+dy^{2}-\left( \ast d\varphi \right) ^{2}, 
\]
where in this case $\varphi \left( x,y\right) $ is the flow potential and $%
\ast $ is the Hodge isomorphism. \ We note that the difference between the
metrics $ds^{\prime 2}$ and $ds^{2}$ corresponds physically to a difference
between a composite metric with noneuclidean part conformally equivalent to
a metric on streamlines, and a composite metric with noneuclidean part
conformally equivalent to a metric on potential lines. \ This correspondence
arises from relating the differential of the \textit{stream function} $\psi $
to the differential of the flow potential $\varphi $ by the equation 
\[
d\psi =c^{2/\left( \gamma _{a}-1\right) }\ast d\varphi . 
\]

\subsection{Cauchy-Riemann equations}

An alternative to considering the functions $u_{1},u_{2}$ to be components
of a 1-form in $\Bbb{R}^{2}$ is to treat them as components of a function in 
$\Bbb{C}$. \ This is a standard approach in which, for example, the
continuity equations in the hodograph plane are associated with a
generalized Cauchy-Riemann operator. \ Among its many advantages, this
approach has the disadvantage of giving special emphasis to dimension 2 and
to the conformal group (or to quasiconformal mappings in the quasilinear
case). \ In fact, the natural invariance group for eqs. (6)-(8) is the
projective group rather than the conformal group, a circumstance which has
some interesting consequences. \ For instance, whereas there are many conic
sections in $\Bbb{R}^{2}$, the unit circle is one of only a few conic
sections in the real projective plane; so the parabolic degeneracy at the
point at infinity in the projective metric corresponds under projective
mappings to a variety of parabolic curves in a euclidean metric (\textit{c.f.%
} [7], Sec. V.86; [3], p. 633).

\section{A boundary-value problem}

The Dirichlet problem for the systems introduced in the preceding section
involves prescribing the value of the 1-form $u_{1}dx+u_{2}dy$ on the
boundary of a domain of $\Bbb{R}^{2}$. \ In the following we consider an
analogue of the Dirichlet problem in which we show the existence of weak
solutions to (6)-(8) which satisfy the boundary condition 
\begin{equation}
u_{1}\frac{dx}{ds}+u_{2}\frac{dy}{ds}=0,
\end{equation}
where $s$ denotes arc length, on the noncharacteristic part of the domain
boundary. \ The proof is based on methods introduced in Ref. 5 for
boundary-value problems in the Chaplygin model.

Denote by $R$ be the region bounded by the rectangle $1/\sqrt{2}<x\leq 1,$ $%
-1/\sqrt{2}<y<1/\sqrt{2}.$ \ Let $C$ be any smooth curve lying entirely in
the interior of $R$ except for two distinct points, which intersect the
characteristic line $x=1$ at $\left( 1,y_{0}\right) $ and $\left(
1,y_{1}\right) ,$ $-1/\sqrt{2}<y_{0}<y_{1}<1/\sqrt{2}$. \ Define $\Omega $
to be the domain bounded by $C\cup \Gamma ,$ where $\Gamma $ is the line
segment $\left( 1,y_{0}\right) \leq \left( x,y\right) \leq \left(
1,y_{1}\right) .$ \ Assume that $dy\leq 0$ on $C$.

The domain $\Omega $ may seem to be rather small and special, but it is not
when the comparison is made to other systems which change type along a conic
section. \ For example, the existence of weak solutions to the Frankl'
problem for the cold plasma model, which changes type along a parabola in $%
\Bbb{R}^{2}$, has been proven only inside a very specific domain contained
within an arbitrarily small circle tangent to the origin [8].

Define $U$ to be the vector space consisting of all pairs of measurable
functions $u=\left( u_{1},u_{2}\right) $ for which the weighted $L^{2}$ norm 
\[
\left\| u\right\| _{\ast }=\left[ \int \int_{\Omega }\left( \left|
2x^{2}-1\right| u_{1}^{2}+\left| 2y^{2}-1\right| u_{2}^{2}\right) dxdy\right]
^{1/2} 
\]
is finite. \ Denote by $W$ the linear space defined by pairs of functions $%
w=\left( w_{1},w_{2}\right) $ having continuous derivatives and satisfying: 
\[
w_{1}dx+w_{2}dy=0 
\]
on $\Gamma $; 
\[
w_{1}=0 
\]
on $C;$%
\[
\int \int_{\Omega }\left[ \left| 2x^{2}-1\right| ^{-1}\left( L^{\ast
}w\right) _{1}^{2}+\left| 2y^{2}-1\right| ^{-1}\left( L^{\ast }w\right)
_{2}^{2}\right] dxdy<\infty . 
\]
Here 
\[
\left( L^{\ast }w\right) _{1}=\left[ \left( 1-x^{2}\right) w_{1}\right]
_{x}-2xyw_{1y}+\left[ \left( 1-y^{2}\right) w_{2}\right] _{y}+2xw_{1}, 
\]
and 
\[
\left( L^{\ast }w\right) _{2}=\left( 1-y^{2}\right) \left(
w_{1y}-w_{2x}\right) +2yw_{1}. 
\]
Define the Hilbert space $H$ to consist of pairs of measurable functions $%
h=\left( h_{1},h_{2}\right) $ for which the norm 
\[
\left\| h\right\| ^{\ast }=\left[ \int \int_{\Omega }\left( \left|
2x^{2}-1\right| ^{-1}h_{1}^{2}+\left| 2y^{2}-1\right| ^{-1}h_{2}^{2}\right)
dxdy\right] ^{1/2} 
\]
is finite.

If the curve $C$ is chosen so that $x$ is bounded below away from the value $%
1/\sqrt{2}$ and $y$ is bounded above and below away from the values $\pm 
\sqrt{1/2},$ then the above weighted inner products can all be replaced by
the $L^{2}$ inner product.

\bigskip

\textbf{Definition}. \ We say that $u$ is a \textit{weak solution} of the
system (6)-(8), (11) in $\Omega $ if $u\in U$ and for every $w\in W,$%
\begin{equation}
-\left( w,g\right) =\left( L^{\ast }w,u\right) ,
\end{equation}
where 
\[
\left( w,g\right) =\int \int_{\Omega }\left( w_{1}g_{1}+w_{2}g_{2}\right)
dxdy. 
\]

\bigskip

The following proposition shows that this notion of weak solution is
well-defined.

\begin{proposition}
Any continuously differentiable weak solution of the boundary-value problem
(6)-(8), (11) with $g\in H$ is a classical solution of the system (6)-(8),
with (11) satisfied on the noncharacteristic curve $C$.
\end{proposition}

\textit{Proof}. \ In the interest of generality, we prove the proposition by
an argument that applies to any smooth domain having a characteristic line
segment on the boundary; we do not use any of the special properties of the
line $x=1$ or of the first and fourth quadrants. 
\[
\left( L^{\ast }w\right) _{1}u_{1}=\left[ \left( 1-x^{2}\right) w_{1}\right]
_{x}u_{1}-2xyw_{1y}u_{1}+\left[ \left( 1-y^{2}\right) w_{2}\right] _{y}u_{1} 
\]
\[
+2xw_{1}u_{1}=\left[ \left( 1-x^{2}\right) w_{1}u_{1}\right] _{x}-\left(
1-x^{2}\right) w_{1}u_{1x} 
\]
\[
-\left[ 2xyw_{1}u_{1}\right] _{y}+2xw_{1}u_{1}+2xyw_{1}u_{1y} 
\]
\[
+\left[ \left( 1-y^{2}\right) w_{2}u_{1}\right] _{y}-\left( 1-y^{2}\right)
w_{2}u_{1y}+2xw_{1}u_{1}, 
\]
and 
\[
\left( L^{\ast }w\right) _{2}u_{2}=\left( 1-y^{2}\right) \left(
w_{1y}-w_{2x}\right) u_{2}+2yw_{1}u_{2} 
\]
\[
=\left[ \left( 1-y^{2}\right) w_{1}u_{2}\right] _{y}+2yw_{1}u_{2}-\left(
1-y^{2}\right) w_{1}u_{2y} 
\]
\[
-\left[ \left( 1-y^{2}\right) w_{2}u_{2}\right] _{x}+\left( 1-y^{2}\right)
w_{2}u_{2x}+w_{1}2yu_{2}, 
\]
Application of Green's Theorem to the derivatives of products yields 
\[
\int \int_{\Omega }\left[ \left( 1-x^{2}\right) w_{1}u_{1}-\left(
1-y^{2}\right) w_{2}u_{2}\right] _{x}dxdy- 
\]
\[
\int \int_{\Omega }\left[ 2xyw_{1}u_{1}-\left( 1-y^{2}\right) \left(
w_{2}u_{1}+w_{1}u_{2}\right) \right] _{y}dxdy= 
\]
\[
\int_{\partial \Omega }\left[ \left( 1-x^{2}\right) w_{1}u_{1}-\left(
1-y^{2}\right) w_{2}u_{2}\right] dy+ 
\]
\[
\int_{\partial \Omega }\left[ 2xyw_{1}u_{1}-\left( 1-y^{2}\right) \left(
w_{2}u_{1}+w_{1}u_{2}\right) \right] dx. 
\]
On the characteristic line segment $\Gamma $ this integral splits into the
sum $I_{1}+I_{2}$, where 
\[
I_{1}=\int_{\Gamma }\left( 1-x^{2}\right) w_{1}u_{1}dy+\left[
2xyw_{1}u_{1}-\left( 1-y^{2}\right) w_{2}u_{1}\right] dx= 
\]
\[
\int_{\Gamma }\left[ 2xyw_{1}u_{1}+\left( 1-x^{2}\right) w_{1}u_{1}\left( 
\frac{dy}{dx}\right) -\left( 1-y^{2}\right) w_{2}u_{1}\right] dx, 
\]
and 
\[
I_{2}=-\int_{\Gamma }\left( 1-y^{2}\right) u_{2}\left(
w_{1}dx+w_{2}dy\right) . 
\]
The integral $I_{2}$ vanishes by the boundary condition for elements of $W,$
from which we also obtain 
\begin{equation}
I_{1}=\int_{\Gamma }\left\{ 2xyw_{1}u_{1}-\left[ \left( 1-x^{2}\right)
\left( \frac{dy}{dx}\right) ^{2}+\left( 1-y^{2}\right) \right]
w_{2}u_{1}\right\} dx.
\end{equation}
On the characteristic curves, 
\[
\left( 1-y^{2}\right) dx^{2}+2xydxdy+\left( 1-x^{2}\right) dy^{2}=0, 
\]
so 
\begin{equation}
\left( 1-x^{2}\right) \frac{dy^{2}}{dx^{2}}=-\left( 1-y^{2}\right) -2xy\frac{%
dy}{dx}.
\end{equation}
Substituting (14) into (13) yields 
\[
I_{1}=\int_{\Gamma }\left\{ 2xyw_{1}u_{1}-\left[ -\left( 1-y^{2}\right) -2xy%
\frac{dy}{dx}+\left( 1-y^{2}\right) \right] w_{2}u_{1}\right\} dx 
\]
\[
=\int_{\Gamma }2xyu_{1}\left( w_{1}+w_{2}\frac{dy}{dx}\right) dx=0. 
\]
Because $w_{1}$ vanishes on $C$, the boundary integral there has the form 
\[
-\int_{C}\left( 1-y^{2}\right) w_{2}\left( u_{1}dx+u_{2}dy\right) . 
\]
We obtain 
\begin{equation}
\left( L^{\ast }w,u\right) =-\left( w,Lu\right) -\int_{C}\left(
1-y^{2}\right) w_{2}\left( u_{1}dx+u_{2}dy\right) ,
\end{equation}
where 
\[
\left( w,Lu\right) = 
\]
\[
\int \int_{\Omega }\left[ \left( 1-x^{2}\right)
w_{1}u_{1x}-2xw_{1}u_{1}-2xyw_{1}u_{1y}+\left( 1-y^{2}\right)
w_{2}u_{1y}-2xw_{1}u_{1}-2yw_{1}u_{2}\right] dxdy 
\]
\[
-\int \int_{\Omega }\left[ 2yw_{1}u_{2}-\left( 1-y^{2}\right)
w_{1}u_{2y}+\left( 1-y^{2}\right) w_{2}u_{2x}\right] dxdy= 
\]
\[
\int \int_{\Omega }\left\{ \left[ \left( 1-x^{2}\right) u_{1}\right]
_{x}-2xyu_{1y}+\left[ \left( 1-y^{2}\right) u_{2}\right]
_{y}-2xu_{1}-2yu_{2}\right\} w_{1}dxdy 
\]
\[
+\int \int_{\Omega }\left( 1-y^{2}\right) \left( u_{1y}-u_{2x}\right)
w_{2}dxdy 
\]
\[
=\int \int_{\Omega }\left[ \left( Lu\right) _{1}w_{1}+\left( Lu\right)
_{2}w_{2}\right] dxdy. 
\]
Combining eqs. (12) and (15) yields 
\[
-\left( w,g\right) =\left( L^{\ast }w,u\right) = 
\]
\[
-\left( w,Lu\right) -\int_{C}\left( 1-y^{2}\right) w_{2}\left(
u_{1}dx+u_{2}dy\right) . 
\]
Because $w$ is arbitrary in $W$, we conclude that (11) is satisfied on $C$
and $Lu=g$, which completes the proof.

In Secs. 4 and 5 we prove:

\begin{theorem}
There exists a weak solution of the boundary-value problem (6)-(8), (11) on $%
\Omega $ for every $g\in H.$
\end{theorem}

\textbf{Remark}. \ Switching the sign of the term $2yu_{2}$ in eq. (7) has
no effect on the proof of Theorem 2.

\section{An \textit{a priori} estimate}

\begin{lemma}
$\exists \,K\in \Bbb{R}^{+}\backepsilon \forall w\in W,$ $K\left\| w\right\|
_{\ast }\leq \left\| L^{\ast }w\right\| ^{\ast }.$
\end{lemma}

\textit{Proof}. \ We use an abbreviated version of the Friedrichs \textit{abc%
} method. \ Fixing a sufficiently differentiable function $a\left(
x,y\right) ,$ consider the $L^{2}$ inner product 
\[
\left( L^{\ast }w,aw\right) = 
\]
\[
\int \int_{\Omega }\left\{ \left[ \left( 1-x^{2}\right) w_{1}\right]
_{x}-2xyw_{1y}+\left[ \left( 1-y^{2}\right) w_{2}\right] _{y}+2xw_{1}\right%
\} aw_{1}dxdy 
\]
\[
+\int \int_{\Omega }\left[ \left( 1-y^{2}\right) \left( w_{1y}-w_{2x}\right)
+2yw_{1}\right] aw_{2}dxdy=\int \int_{\Omega _{m}}\sum_{i=1}^{7}\tau
_{i}dxdy, 
\]
where 
\begin{equation}
\tau _{1}=\frac{1}{2}\left[ \left( 1-x^{2}\right) aw_{1}^{2}\right] _{x}-%
\left[ \frac{1}{2}\left( 1-x^{2}\right) a_{x}+ax\right] w_{1}^{2};
\end{equation}
\begin{equation}
\tau _{2}=-\left( xyaw_{1}^{2}\right) _{y}+\left( ax+xya_{y}\right)
w_{1}^{2};
\end{equation}
\begin{equation}
\tau _{3}=\left[ \left( 1-y^{2}\right) aw_{1}w_{2}\right] _{y}-\left(
1-y^{2}\right) a_{y}w_{2}w_{1}-\left( 1-y^{2}\right) aw_{2}w_{1y};
\end{equation}
\begin{equation}
\tau _{4}=2xaw_{1}^{2};\;\tau _{5}=\left( 1-y^{2}\right) aw_{1y}w_{2};
\end{equation}
\begin{equation}
\tau _{6}=-\frac{1}{2}\left[ \left( 1-y^{2}\right) aw_{2}^{2}\right] _{x}+%
\frac{1}{2}\left( 1-y^{2}\right) a_{x}w_{2}^{2};
\end{equation}
\begin{equation}
\;\tau _{7}=2yaw_{1}w_{2}.
\end{equation}

We ignore for a moment derivatives of products, as these will be integrated
and become boundary terms. \ The coefficients of $w_{1y}w_{2}$ sum to zero
in (18) and (19). \ Denoting the coefficients of $w_{1}^{2}$ off the
boundary by $\alpha ,$ those of $w_{2}^{2}$ by $\gamma $ and those of $%
w_{1}w_{2}$ by $2\beta $ and choosing $a=x^{2},$ we obtain 
\[
\alpha =x\left( 3x^{2}-1\right) ; 
\]
\[
\gamma =x\left( 1-y^{2}\right) ; 
\]
\[
2\beta =2yx^{2}. 
\]
The region $R$ is defined so that the discriminant 
\[
\alpha \gamma -\beta ^{2}=x^{2}\left[ y^{2}\left( 1-4x^{2}\right) +3x^{2}-1%
\right] 
\]
is positive on $\Omega $. \ Thus we have the estimate 
\[
2\beta w_{1}w_{2}\geq -2\left| \beta \right| \left| w_{1}\right| \left|
w_{2}\right| >-2\sqrt{\alpha }\left| w_{1}\right| \sqrt{\gamma }\left|
w_{2}\right| \geq -\alpha w_{1}^{2}-\gamma w_{2}^{2}. 
\]
This shows that the inequality of the lemma is satisfied in $\Omega ,$ but
the resulting constant depends on the choice of the curve $C$. \ Rather, we
prefer to obtain the explicit estimate 
\[
2\beta w_{1}w_{2}\geq -2x\left| xw_{1}\right| \left| yw_{2}\right| \geq
-\left( x^{3}w_{1}^{2}+xy^{2}w_{2}^{2}\right) . 
\]

Applying Green's Theorem to derivatives of products in $\left( Lw,aw\right) $
results in a boundary integral of the form 
\[
\int_{\partial \Omega }\frac{x^{2}}{2}\left[ \left( 1-x^{2}\right)
w_{1}^{2}-\left( 1-y^{2}\right) w_{2}^{2}\right] dy+ 
\]
\[
\int_{\partial \Omega }x^{2}\left[ xyw_{1}^{2}-\left( 1-y^{2}\right)
w_{1}w_{2}\right] dx. 
\]
The definition of $W$ implies that on $\Gamma ,$ 
\[
-\left( 1-y^{2}\right) w_{1}w_{2}dx=\left( 1-y^{2}\right) w_{2}^{2}dy, 
\]
so the boundary integral on $\Gamma $ reduces to 
\[
\int_{\Gamma }\frac{x^{2}}{2}\left[ \left( 1-x^{2}\right) w_{1}^{2}+\left(
1-y^{2}\right) w_{2}^{2}\right] dy+x^{3}yw_{1}^{2}dx 
\]
\[
=\int_{\Gamma }\frac{x^{2}}{2}\left( 1-y^{2}\right) w_{2}^{2}dy=0. 
\]
Because $w_{1}$ vanishes on $C,$ the remaining boundary integral is of the
form 
\[
-\int_{C}\frac{x^{2}}{2}\left( 1-y^{2}\right) w_{2}^{2}dy, 
\]
which is nonnegative under the given orientation by the hypotheses on $C$.

We find that on $\Omega ,$ 
\begin{equation}
\left( L^{\ast }w,aw\right) \geq \frac{1}{\sqrt{2}}\int \int_{\Omega }\left(
\left| 2x^{2}-1\right| w_{1}^{2}+\left| 2y^{2}-1\right| w_{2}^{2}\right)
dxdy.
\end{equation}

It remains to estimate $\left( L^{\ast }w,aw\right) $ from above. \ We have
for any positive constant $\lambda $, 
\[
\left( L^{\ast }w,aw\right) \leq \int \int_{\Omega }\left| \sqrt{2x^{2}-1}%
w_{1}\right| \left| \left( \sqrt{2x^{2}-1}\right) ^{-1}\left( L^{\ast
}w\right) _{1}\right| dxdy 
\]
\[
+\int \int_{\Omega }\left| \sqrt{\left| 2y^{2}-1\right| }w_{2}\right| \left|
\left( \sqrt{\left| 2y^{2}-1\right| }\right) ^{-1}\left( L^{\ast }w\right)
_{2}\right| dxdy 
\]
\begin{equation}
\leq \frac{1}{\lambda }\left\| L^{\ast }w\right\| ^{\ast 2}+\lambda \left\|
w\right\| _{\ast }^{2}.
\end{equation}
Choosing $\lambda <1/\sqrt{2},$ inequalities (22) and (23) imply the
assertion of Lemma 3 with $K=\sqrt{\left[ \left( 1/\sqrt{2}\right) -\lambda %
\right] \lambda }.$

\section{Existence}

The proof of existence is straightforward, given the \textit{a priori}
estimates of the preceding section. \ We briefly outline the argument,
following Ref. 5.

Define the scaled 1-forms 
\[
\widetilde{w}=\sqrt{2x^{2}-1}w_{1}dx+\sqrt{\left| 2y^{2}-1\right| }w_{2}dy 
\]
and 
\[
\widetilde{g}=\frac{1}{\sqrt{2x^{2}-1}}g_{1}dx+\frac{1}{\sqrt{\left|
2y^{2}-1\right| }}g_{2}dy. 
\]
Arguing as in (23), but applying the Schwartz inequality in place of Young's
inequality, we have 
\[
\left| \left( w,g\right) \right| =\left| \left( \widetilde{w},\widetilde{g}%
\right) \right| \leq \left\| \widetilde{w}\right\| _{2}\left\| \widetilde{g}%
\right\| _{2}, 
\]
where $\left\| \;\right\| _{2}$ is the (unweighted) $L^{2}$ norm. \ The
extreme left- and right-hand sides of this inequality can be written 
\[
\left| \left( w,g\right) \right| \leq \left\| w\right\| _{\ast }\left\|
g\right\| ^{\ast }\leq 
\]
\[
K^{-1}\left\| L^{\ast }w\right\| ^{\ast }\left\| g\right\| ^{\ast }\leq 
\widetilde{K}\left( g\right) \left\| L^{\ast }w\right\| ^{\ast }, 
\]
using Lemma 3. \ Thus the functional $\xi $ defined for fixed $g$ and all $%
w\in W$ by the formula 
\[
\xi \left( L^{\ast }w\right) =-\left( w,g\right) 
\]
can be extended to a bounded linear functional on $H$. \ The Riesz
Representation Theorem implies that $\forall w\in W$ there is an $h\in H$
for which 
\[
\xi \left( L^{\ast }w\right) =\left( L^{\ast }w,h\right) ^{\ast }. 
\]
Defining $u=\left( u_{1},u_{2}\right) $ so that 
\[
u_{1}=-\left( 2x^{2}-1\right) ^{-1}h_{1} 
\]
and 
\[
u_{2}=-\left| 2y^{2}-1\right| ^{-1}h_{2}, 
\]
we have $u\in U$ as $h\in H;$ that is, 
\[
\int \int_{\Omega }\left[ \left( 2x^{2}-1\right) u_{1}^{2}+\left|
2y^{2}-1\right| u_{2}^{2}\right] dxdy= 
\]
\[
\int \int_{\Omega }\left[ \left( 2x^{2}-1\right) ^{-1}h_{1}^{2}+\left|
2y^{2}-1\right| ^{-1}h_{2}^{2}\right] dxdy<\infty . 
\]

We conclude that 
\[
-\left( w,g\right) =\xi \left( L^{\ast }w\right) =\left( L^{\ast }w,h\right)
^{\ast }= 
\]
\[
\int \int_{\Omega }\left[ \left( 2x^{2}-1\right) ^{-1}\left( L^{\ast
}w\right) _{1}h_{1}+\left| 2y^{2}-1\right| ^{-1}\left( L^{\ast }w\right)
_{2}h_{2}\right] dxdy= 
\]
\[
-\int \int_{\Omega }\left[ \left( 2x^{2}-1\right) ^{-1}\left( L^{\ast
}w\right) _{1}\left( 2x^{2}-1\right) u_{1}+\left| 2y^{2}-1\right|
^{-1}\left( L^{\ast }w\right) _{2}\left| 2y^{2}-1\right| u_{2}\right] dxdy 
\]
\[
=\left( L^{\ast }w,u\right) . 
\]

Comparing the extreme left-hand side of this expression with its extreme
right-hand side completes the proof of Theorem 2.

\section{Modifications of the problem}

The lower-order terms of equations of mixed type are frequently modified in
order to simplify the analysis [see, for example, eqs. (7) and (23) of Ref.
6 or eqs. (1.11) and (2.1) of Ref. 8]. \ In addition to solving the system
(6)-(8), we can also prove the existence of weak solutions to a systems
which differ from (6)-(8) only in the form of their lower-order terms. \
Among many possible examples, we choose two obvious ones.

\subsection{A different distribution of the lower-order terms}

We can replace (6)-(8) with a system having the more symmetric form 
\[
\widetilde{L}u=g, 
\]
where 
\[
\widetilde{L}=\left( \widetilde{L}_{1},\widetilde{L}_{2}\right) ,\;g=\left(
g_{1},g_{2}\right) , 
\]
\[
\left( \widetilde{L}u\right) _{1}=\left[ \left( 1-x^{2}\right) u_{1}\right]
_{x}-2xyu_{1y}+\left[ \left( 1-y^{2}\right) u_{2}\right] _{y}-2xu_{1}, 
\]
and 
\[
\left( \widetilde{L}u\right) _{2}=\left( 1-y^{2}\right) \left(
u_{1y}-u_{2x}\right) +2yu_{2}. 
\]
\ In the special case $g_{1}=g_{2}=0$ both this system and eqs. (6)-(8)
satisfy the equation 
\[
\left[ \left( 1-x^{2}\right) u_{1}\right] _{x}-2xyu_{1y}+\left[ \left(
1-y^{2}\right) u_{2}\right] _{y}-2\left( xu_{1}+yu_{2}\right) =\left(
1-y^{2}\right) \left( u_{1y}-u_{2x}\right) , 
\]
although the equated quantities differ in the different systems. \ The
analysis of the modified system is a little simpler and the conditions on
the noncharacteristic part of the boundary considerably more lenient. \
However, the proof for this system does not apply in an obvious way to a
domain lying in two contiguous quadrants.

Denote by $\Omega _{m}$ the region bounded by the characteristic line
tangent to the unit disc at the point $\left( 1,0\right) $ and a smooth
curve $C_{m}$ which intersects that line at exactly two points on the line
segment $\Gamma $ given by the interval $\left( 1,-1\right) <\left(
1,y\right) <\left( 1,0\right) .$ \ Assume that $C_{m}$ is bounded on the
left by the line $x=0$, on the right by $\Gamma ,$ below by the line $y=-1,$
and above by the $x$-axis. \ Orient $\partial \Omega _{m}$ in the
counterclockwise direction. \ We assume that, with this orientation, the
line element $dy$ is nonpositive on $C_{m}$. \ (Small modifications of the
problem will define an analogous boundary-value problem in the second
quadrant, a fact which is reflected below in our notation for the spaces $%
U_{m}$, $W_{m}$, and $H_{m}$.)

Denote by $U_{m}$ the vector space consisting of all pairs of measurable
functions $u=\left( u_{1},u_{2}\right) $ for which the weighted $L^{2}$ norm 
\[
\left\| u\right\| _{m\ast }=\left\{ \int \int_{\Omega _{m}}\left( \left|
x\right| u_{1}^{2}+\left| y\right| u_{2}^{2}\right) dxdy\right\} ^{1/2} 
\]
is finite. \ This norm is induced by the weighted inner product 
\[
\left( u,w\right) _{m\ast }=\int \int_{\Omega _{m}}\left( \left| x\right|
u_{1}w_{1}+\left| y\right| u_{2}w_{2}\right) dxdy. 
\]

Denote by $W_{m}$ the linear space defined by pairs of functions $w=\left(
w_{1},w_{2}\right) $ having continuous derivatives and satisfying: 
\[
w_{1}dx+w_{2}dy=0 
\]
on $\Gamma ;$%
\[
w_{1}=0 
\]
on $C_{m}$; 
\[
\int \int_{\Omega _{m}}\left[ \left| x\right| ^{-1}\left( \widetilde{L}%
^{\ast }w\right) _{1}^{2}+\left| y\right| ^{-1}\left( \widetilde{L}^{\ast
}w\right) _{2}^{2}\right] dxdy<\infty . 
\]
Here 
\[
\left( \widetilde{L}^{\ast }w\right) _{1}=\left[ \left( 1-x^{2}\right) w_{1}%
\right] _{x}-2xyw_{1y}+\left[ \left( 1-y^{2}\right) w_{2}\right]
_{y}+2xw_{1}, 
\]
and 
\[
\left( \widetilde{L}^{\ast }w\right) _{2}=\left( 1-y^{2}\right) \left(
w_{1y}-w_{2x}\right) -2yw_{2}. 
\]
The space $W_{m}$ is contained in the Hilbert space $H_{m}$ consisting of
pairs of measurable functions $h=\left( h_{1},h_{2}\right) $ for which the
norm 
\[
\left\| h\right\| _{m}^{\ast }=\left\{ \int \int_{\Omega _{m}}\left( \left|
x\right| ^{-1}h_{1}^{2}+\left| y\right| ^{-1}h_{2}^{2}\right) dxdy\right\}
^{1/2} 
\]
is finite.

To prove the analogue of Lemma 3 for this system under an analogous boundary
condition we estimate the $L^{2}$ inner product $\left( w,\widetilde{L}%
^{\ast }w\right) $ as in (16)-(21). \ We obtain a result analogous to (22)
with

\[
\alpha =2x,\;\gamma =-2y, 
\]
and 
\[
2\beta =0. 
\]
Arguing as in (23) with $\lambda <2,$ we find that

\[
\exists \,K_{m}\in \Bbb{R}^{+}\backepsilon \forall w\in W_{m},K_{m}\left\|
w\right\| _{m\ast }\leq \left\| \widetilde{L}^{\ast }w\right\| _{m}^{\ast } 
\]
with $K_{m}=\sqrt{\left( 2-\lambda \right) \lambda }.$ \ The remainder of
the existence proof proceeds as in the case of eqs. (6)-(8).

If the term $2yu_{2}$ in the component $\left( \widetilde{L}u\right) _{2}$
is multiplied by $-1$, then the domain of the solution switches from the
fourth quadrant to the first quadrant. \ If the term $2xu_{1}$ in the
component $\left( \widetilde{L}u\right) _{1}$\ is multiplied by $-1$, then
the domain switches from the fourth quadrant to the third quadrant. \ If
both lower-order terms are multiplied by $-1$, then the domain switches from
the fourth quadrant to the second quadrant. \ In the last two cases, $\Gamma 
$ lies along the line $x=-1$.

\subsection{Neglected lower-order terms}

Finally, we consider a form of the system (6)-(8) in which no terms of order
zero appear. \ This system consists of equations having the form 
\[
L_{o}u=g, 
\]
where 
\[
L_{o}=\left( L_{o1},L_{o2}\right) ,\;g=\left( g_{1},g_{2}\right) , 
\]
\[
\left( L_{o}u\right) _{1}=\left[ \left( 1-x^{2}\right) u_{1}\right]
_{x}-2xyu_{1y}+\left[ \left( 1-y^{2}\right) u_{2}\right] _{y}, 
\]
and 
\[
\left( L_{o}u\right) _{2}=\left( 1-y^{2}\right) \left( u_{1y}-u_{2x}\right)
. 
\]
In this case the boundary-value problem is simplified somewhat by the fact
that $L_{o}=L_{o}^{\ast }.$ \ For example, Lemma 3 implies the uniqueness in 
$W$ of weak solutions, which are defined by direct analogy to the other two
cases.

To prove the existence of weak solutions to the system $L_{o}u=g,$ we fix
positive numbers $\delta <<1/2$ and $\varepsilon <<1/2$ and denote by $R_{o}$
the rectangle 
\[
\frac{1}{\sqrt{2}}<x\leq 1,\;\frac{1}{\sqrt{2-\delta }}<y\leq \sqrt{%
1-\varepsilon }. 
\]
\ Let $C_{o}$ be a smooth curve lying in the interior of $R_{o}$ with the
exception of two distinct points, $\left( 1,y_{0}\right) $ and $\left(
1,y_{1}\right) ,$ $1/\sqrt{2-\delta }<y_{0}<y_{1}\leq \sqrt{1-\varepsilon },$
at which the curve intersects the characteristic line $x=1$. \ Define $%
\Gamma $ to be the line segment $\left( 1,y_{0}\right) $ $\leq \left(
x,y\right) \leq \left( 1,y_{1}\right) $ and $\Omega _{o}$ to be the domain
having boundary $C_{o}\cup \Gamma .$ \ In this case we can take the
associated Hilbert spaces, $U_{o}$ and $H_{o},$ to be $L^{2}$, bearing in
mind that our estimates will depend in a predictable way on the sizes of $%
\varepsilon $ and $\delta .$ \ As in the preceding cases, we place the
boundary condition (11) on the noncharacteristic part of the boundary.

In order to prove the analogue of Lemma 3 for this system, we estimate the $%
L^{2}$ inner product 
\[
\left( L_{o}w,xyw\right) =\int \int_{\Omega _{o}}\sum_{i=1}^{5}\tau _{i}dxdy 
\]
where 
\[
\tau _{1}=\frac{1}{2}\left[ \left( 1-x^{2}\right) xyw_{1}^{2}\right]
_{x}-x^{2}yw_{1}^{2}-\frac{1}{2}\left( 1-x^{2}\right) yw_{1}^{2}; 
\]
\[
\tau _{2}=-\left[ x^{2}y^{2}w_{1}^{2}\right]
_{y}+x^{2}yw_{1}^{2}+x^{2}yw_{1}^{2}; 
\]
\[
\tau _{3}=\left[ \left( 1-y^{2}\right) xyw_{1}w_{2}\right] _{y}-\left(
1-y^{2}\right) xw_{1}w_{2}-\left( 1-y^{2}\right) xyw_{1y}w_{2}; 
\]
\[
\tau _{4}=\left( 1-y^{2}\right) w_{1y}xyw_{2}; 
\]
\[
\tau _{5}=-\frac{1}{2}\left[ \left( 1-y^{2}\right) xyw_{2}^{2}\right] _{x}+%
\frac{1}{2}\left( 1-y^{2}\right) yw_{2}^{2}. 
\]
We have 
\[
\alpha =\frac{y}{2}\left( 3x^{2}-1\right) ,\;\gamma =\frac{y}{2}\left(
1-y^{2}\right) , 
\]
and 
\[
2\beta =-\left( 1-y^{2}\right) x, 
\]
yielding 
\[
\alpha \gamma -\beta ^{2}=\frac{y^{2}\left( 1-y^{2}\right) }{4}\left( 3x^{2}-%
\frac{1-y^{2}}{y^{2}}x^{2}-1\right) . 
\]
Because $R_{o}$ is constructed so that 
\begin{equation}
\frac{1-y^{2}}{y^{2}}<1-\delta <1,
\end{equation}
it is sufficient to show that 
\[
2x^{2}-1>0, 
\]
which also follows from the construction of $R_{o}$. \ Now 
\begin{equation}
-2\beta w_{1}w_{2}\geq -\frac{\sqrt{1-y^{2}}}{2}\left[ x^{2}w_{1}^{2}+\left(
1-y^{2}\right) w_{2}^{2}\right] .
\end{equation}
Taking square roots in (24), and applying the result to (25) yields 
\[
-2\beta w_{1}w_{2}>-\frac{\sqrt{1-\delta }}{2}\left[ yx^{2}w_{1}^{2}+y\left(
1-y^{2}\right) w_{2}^{2}\right] . 
\]
Thus we have 
\[
\left( L_{o}w,xyw\right) \geq \frac{\varepsilon \left( 1-\sqrt{1-\delta }%
\right) }{2\sqrt{2-\delta }}\left\| w\right\| _{2}^{2}. 
\]
The remainder of the existence proof is exactly analogous to the arguments
for the preceding cases.

\bigskip

\textbf{Acknowledgment}. \ I am indebted to an anonymous referee for helpful
criticism of an earlier draft of this paper.

\bigskip

\textbf{References}

[1] Beltrami, E., Teoria fondamentale degli spazii di curvatura costante,
Annali di Matematica Pura ed Applicata, ser. 2, \textbf{2}, 232-255 (1868).

[2] Bers, L., Mathematical Aspects of Subsonic and Transonic Gas Dynamics,
New York: Wiley, 1958.

[3] Hua, L. K. \ In: Proceedings of the 1980 Beijing Symposium on
Differential Geometry and Differential Equations. S. S. Chern and Wu
Wen-ts\"{u}n (eds.), pp. 627-654, New York: Gordon and Breach, 1982.

[4] Ji Xin-hua and Chen De-quan, Tricomi's problems of non-homogeneous
equation of mixed type in real projective plane, in: \textit{Proceedings of
the 1980 Beijing Symposium on Differential Geometry and Differential
Equations,} (S. S. Chern and Wu Wen-ts\"{u}n, eds.), pp. 1257-1271, Gordon
and Breach, New York, 1982.

[5] Morawetz, C. S., A weak solution for a system of equations of
elliptic-hyperbolic type, Commun. Pure Appl. Math. \textbf{11}, 315-331
(1958).

[6] Morawetz, C. S., Stevens, D. C., and Weitzner, H., A numerical
experiment on a second-order partial differential equation of mixed type,
Commun. Pure Appl. Math. \textbf{44}, 1091-1106 (1991).

[7] Veblen, O., and Young, J. W., Projective Geometry, Vol. II, Boston: Ginn
and Co., 1918.

[8] Yamamoto, Y., Existence and uniqueness of a generalized solution for a
system of equations of mixed type. Ph.D. Dissertation, Polytechnic
University of New York, 1994.

[9] Y. Yang, Classical solutions in the Born-Infeld theory, Proc. R. Soc.
Lond. Ser. A \textbf{456} (2000), no. 1995, 615-640.

\end{document}